\documentstyle [12pt]{article}
\newtheorem{theorem}{Theorem}[section]
\newtheorem{prop}[theorem]{Theorem}

\newtheorem{coroll}[theorem]{Corollary}

\newtheorem{lemm}[theorem]{Lemma}

\newtheorem{anoprop}[theorem]{Proposition}

\newenvironment{mproof}{\begin{trivlist}\item[]{\em
Proof: }}{\hfill$\Box$\end{trivlist}}

\newtheorem{eg}{\rm\sl \uppercase{Example}}[section]

\def \IR{\hbox{{\rm I}\kern-.2em\hbox{{\rm R}}}}
\def \iR{\hbox{{\sevenrm I\kern-.2em\hbox{\sevenrm R}}}}
\def \IN{\hbox{{\rm I}\kern-.2em\hbox{\rm N}}}
\def \IC{\hbox{{\rm I}\kern-.6em\hbox{\bf C}}}
\def \IQ{\hbox{{\rm I}\kern-.6em\hbox{\bf Q}}}
\def \ZZ{\hbox{{\rm Z}\kern-.4em\hbox{\rm Z}}}

\begin{document}
\begin{titlepage}

\vspace*{1cm}
\begin{center}
\Large
{\bf ON THE OUTER AUTOMORPHISM GROUPS OF TRIANGULAR ALTERNATION LIMIT
ALGEBRAS					     }
\end{center}
\vspace*{0.25cm}

\begin{center}
S. C. Power	  \\
Department of Mathematics  \\
University of Lancaster   \\
Lancaster LA1 4YF\\
England.\\
\vspace*{0.50cm}
\end{center}

\bf ABSTRACT

\rm Let $A$
denote the alternation limit
algebra, studied by Hopenwasser and Power,
and by Poon, which is the closed direct limit of upper triangular
matrix algebras determined by refinement embeddings of multiplicity
$r_k$ and standard embeddings of multiplicity $s_k$. It is shown that
the quotient of the isometric automorphism group by the
approximately inner automorphisms is the abelian group
$ \ZZ ^d$ where $d$ is the number of
primes
that are divisors of infinitely many terms of each of the sequences
$(r_k)$ and $(s_k)$. This group is also the group of
automorphisms of the fundamental relation of $A$.

\end{titlepage}

\section{Introduction}

In Hopenwasser and Power [HP] and in Poon [Po] the alternation limit algebras
described below were classified. In this
note we determine the quotient group $Out_{isom}A = Aut_{isom}A/I(A)$ for these
algebras
where $Aut_{isom}A$ is the group of isometric algebra automorphisms and $I(A)$
is the normal subgroup of $AutA$ of approximately inner automorphisms.	An
automorphism $\alpha$ is said to be approximately inner if there exists a
sequence $(b_k)$ of invertible elements such that $\alpha(a) = \lim_k b_k
ab^{-1}_k$ for all $a$ in $A$.

Let $(r_k)$, $(s_k)$ be sequences of positive integers. Write $T(r_k,s_k)$ for
the Banach algebra limit of the system
\[
\IC \to T_{r_1} \to T_{r_1s_1} \to T_{r_1s_1r_2} \to \dots ,
\]
where $T_n$ is the algebra of upper triangular $n\times n$ complex matrices and
where the embeddings are unital and are alternately of refinement type $(\rho
(a) = (a_{ij}1_t)$, with $1_t$ the $t\times t$ identity and of standard type
$(\sigma (a) = a \oplus \dots \oplus a$, $t$ times).

{\bf Theorem 1} $Out_{isom}(T(r_k,s_k)) = \ZZ ^d$ where $d$ is the number of
primes
$p$  that are divisors of infinitely many terms of each of the sequences
$(r_k)$ and $(s_k)$.  (If $d = \infty$ interpret $\ZZ ^d$ as the countably
generated free abelian group.)

The proof uses the methods of [HP]. A major step is to characterise
the automorphism group of the fundamental relation, or semigroupoid, which is
associated with an alternation algebra. This order-topological result is of
independent interest and is stated and proved separately below.

Let $r$ and $s$ be the generalised integers $r_1r_2 \dots,$ and $s_1s_2\dots$
respectively and suppose that $p$ is a prime satisfying the condition in the
statement of the theorem.  Then $p^\infty$ divides $r$ and $s$. Thus we can
arrange new formal products $r = t_1t_2\dots$, $s = u_1u_2 \dots$, with $t_k =
u_k = p$ for all odd $k$.  As noted in [HP], because of the commutation of
refinement and standard embeddings, we can easily display a commuting zig zag
diagram to show that $T(r_k,s_k)$ and $T(t_k,u_k)$ are isometrically
isomorphic.
However, with the new formal product we can construct one of the generators of
$Out_{isom}A$. Consider the automorphism $\alpha$ determined by the following
commuting diagram where the matrix algebras are omitted for notational economy.
\begin{center}
\setlength{\unitlength}{0.0125in}%
\begin{picture}(445,119)(55,700)
\thicklines
\put(130,800){\vector( 1, 0){ 40}}
\put(180,800){\vector( 1, 0){ 40}}
\put(230,800){\vector( 1, 0){ 40}}
\put(280,800){\vector( 1, 0){ 40}}
\put(330,800){\vector( 1, 0){ 40}}
\put(380,800){\vector( 1, 0){ 40}}
\put( 80,720){\vector( 1, 0){ 40}}
\put(130,720){\vector( 1, 0){ 40}}
\put(180,720){\vector( 1, 0){ 40}}
\put(230,720){\vector( 1, 0){ 40}}
\put(280,720){\vector( 1, 0){ 40}}
\put(330,720){\vector( 1, 0){ 40}}
\put(380,720){\vector( 1, 0){ 40}}
\put(140,810){\makebox(0,0)[lb]{\raisebox{0pt}[0pt][0pt]{$\sigma_{u_1}$}}}
\put(190,810){\makebox(0,0)[lb]{\raisebox{0pt}[0pt][0pt]{$\rho_{t_2}$}}}
\put(240,810){\makebox(0,0)[lb]{\raisebox{0pt}[0pt][0pt]{$\sigma_{u_2}$}}}
\put(290,810){\makebox(0,0)[lb]{\raisebox{0pt}[0pt][0pt]{$\rho_{t_3}$}}}
\put(340,810){\makebox(0,0)[lb]{\raisebox{0pt}[0pt][0pt]{$\sigma_{u_3}$}}}
\put( 90,810){\makebox(0,0)[lb]{\raisebox{0pt}[0pt][0pt]{$\rho_{t_1}$}}}
\put(140,700){\makebox(0,0)[lb]{\raisebox{0pt}[0pt][0pt]{$\sigma_{u_1}$}}}
\put( 80,800){\vector( 1, 0){ 40}}
\put(190,700){\makebox(0,0)[lb]{\raisebox{0pt}[0pt][0pt]{$\rho_{t_2}$}}}
\put(355,755){\makebox(0,0)[lb]{\raisebox{0pt}[0pt][0pt]{$\rho_{t_3}$}}}
\put(240,700){\makebox(0,0)[lb]{\raisebox{0pt}[0pt][0pt]{$\sigma_{u_2}$}}}
\put(290,700){\makebox(0,0)[lb]{\raisebox{0pt}[0pt][0pt]{$\rho_{t_3}$}}}
\put(340,700){\makebox(0,0)[lb]{\raisebox{0pt}[0pt][0pt]{$\sigma_{u_3}$}}}
\put( 90,700){\makebox(0,0)[lb]{\raisebox{0pt}[0pt][0pt]{$\rho_{t_1}$}}}
\put( 80,780){\vector( 1,-1){ 35}}
\put(280,780){\vector( 1,-1){ 35}}
\put(130,745){\vector( 1, 1){ 35}}
\put(330,745){\vector( 1, 1){ 35}}
\put(500,780){\vector( 0, -1){ 40}}
\put(440,800){\makebox(0,0)[lb]{\raisebox{0pt}[0pt][0pt]{\twlrm . . . }}}
\put(440,720){\makebox(0,0)[lb]{\raisebox{0pt}[0pt][0pt]{\twlrm . . .}}}
\put(475,715){\makebox(0,0)[lb]{\raisebox{0pt}[0pt][0pt]{$T(t_k,u_k)$}}}
\put(475,795){\makebox(0,0)[lb]{\raisebox{0pt}[0pt][0pt]{$T(t_k,u_k)$}}}
\put( 55,795){\makebox(0,0)[lb]{\raisebox{0pt}[0pt][0pt]{$\IC$}}}
\put( 55,715){\makebox(0,0)[lb]{\raisebox{0pt}[0pt][0pt]{$\IC$}}}
\put( 80,755){\makebox(0,0)[lb]{\raisebox{0pt}[0pt][0pt]{$\sigma_{t_1}$}}}
\put(155,755){\makebox(0,0)[lb]{\raisebox{0pt}[0pt][0pt]{$\rho_{t_1}$}}}
\put(275,755){\makebox(0,0)[lb]{\raisebox{0pt}[0pt][0pt]{$\sigma_{t_3}$}}}
\end{picture}
\end{center}
It will be shown below that $\alpha$ provides a nonzero coset and
that the totality of such cosets provides a generating set for the
isometric outer automorphism group.

\section{Proof of Theorem 1}

Let $X$, or $X(r_k,s_k)$, be the Cantor space
\[
X = \Pi^{-1}_{k=-\infty}\{1,\dots,s_{-k}\} \times \Pi^\infty_{k=1}
\{1, \dots, r_k\},
\]
where we have fixed the sequences $(r_k)$ and $(s_k)$. Define the equivalence
relation $\tilde{R}$ on $X$ to consist of the pairs $(x,y)$ of points $x =
(x_k), y = (y_k)$ in $X$ with $x_k = y_k$ for all large enough and small enough
$k$.  $\tilde{R}$ carries a natural locally compact Hausdorff topology (giving
it
the structure of an approximately finite groupoid).  Write $R$, or
$R(r_k,s_k)$,
for the antisymmetric topologised subrelation of $\tilde{R}$ consisting of
pairs
$(x,y)$ in $R$ for which $x$ preceeds $y$ in the lexicographic order. Thus
$(x,y) \in R$ if and only if $(x,y) \in \tilde{R}$ and, either $x = y$, or for
the smallest $k$ for which $x_k \neq y_k$ we have $x_k < y_k$.

An automorphism of $R(r_k,s_k)$ is a binary relation isomorphism (implemented
by
a bijection $\alpha$ of the underlying space $X$), which is
a homeomorphism for the (relative groupoid) topology of $R(r_k,s_k)$.
Necessarily $\alpha$ is a homeomorphism of $X$.

{\bf Theorem 2} The group of automorphisms of the topological
binary relation
$R(r_k,s_k)$ is $\ZZ ^d$ where $d$  is the number of primes which divide
infinitely many terms of each of the sequences $(r_k)$ and $(s_k)$.

\begin{mproof}	Let $\overline{{\cal O} (x)}$ denote the closure of the
$R$-orbit of the
point $x$ in $X$.  Here ${\cal O}(x) = \{y : (y,x) \in R\}$. Recall from [HP]
that the pair of points $x, x^+$ is called a {\em gap pair} if $x^+ \not\in
\overline{{\cal O}(x)}$ and
\[
\overline{{\cal O} (x^+)} = \overline{{\cal O}(x)} \cup \{x\}.
\]
Furthermore $x,x^+$ is a gap pair if and only if \\
1) there exists $n$ such that $x_m = 1$ for all $m \leq n$,\\
2) there exists $p$ such that $x_q = r_q$ for all $q \geq p$.

 Also if
$p$ is the smallest integer for which 2) holds (with $r_p = s_{-p}$ if $p$ is
negative), then $x^+$ is given by \[
(x^+)_j = \left\{
\begin{array}{l}
x_j \mbox { if } j < p - 1\\
x_{p-1}+1 \mbox { if } j = p - 1\\
1 \mbox { if } j \geq p\\
\end{array}
\right.
\]
The usefulness of this for our purpose is that an automorphism $\alpha$ of $R$
necessarily maps gap pairs to gap pairs and so the coordinate description of
these pairs leads ultimately to a coordinate description of $\alpha$.

Let $\alpha$ be an automorphism of $R$. Consider the (left) gap point $x_\ast =
(\dots, 1, 1, \hat{1},r_1,r_2,\dots)$ where $\hat{1}$ indicates the coordinate
position for $s_1$. Then $\alpha(x_\ast)$ is necessarily a (left) gap point,
thus \[
\alpha(x_\ast) = (\dots, 1, 1, z_{-t+1},z_{-t},\dots, z_{t-1},r_t,r_{t+1}\dots
)
\]
for some positive integer $t$.	We have
\[
\overline{{\cal O}(x_\ast)} = \{x = (\dots, 1, \hat{1},x_1,x_2,\dots) : x_k
\leq
r_k \mbox { for all } k\},
\]
\[
\overline{{\cal O}(\alpha(x_\ast))} = \{y = (\dots 1, w^\prime, y_t,
y_{t+1},\dots)\},
\]
where $y_k \leq r_k$ for all $k \geq t$ and where $w^\prime$ is any word of
length $2t-2$ which preceeds (or is equal to) the word $w = z_{-t+1},z_{-t},
\dots, z_{t-1}$ in the lexicographic order.  Restating this, we have natural
homeomorphisms
\[
\overline{{\cal O}(x_\ast)} \approx \Pi^\infty_{k=1} \{1,\dots,r_k\}
\]
\[
\overline{{\cal O}(\alpha(x_\ast))} \approx \{1, \dots n\} \times
\Pi^\infty_{k=t} \{1, \dots, r_t\}
\]
where $n$ is the number of words $w^\prime$. Moreoever, these identifying
homeomorphisms induce isomorphisms between the restrictions $R|\overline{{\cal
O}(x_\ast)}$ and $R|\overline{{\cal O}(\alpha(x_\ast))}$ and the unilateral
relations $R_1$ and $R_2$, respectively, where $R_1 = R(r_k,u_k)$, with $u_k =
1$ for all $k$, and $R_2 = R(r^\prime_k, u_k)$, with $u_k$ as before, $r^\prime
_1 = n$, and $r^\prime _k = r_{k+t-2}$ for $k = 2,3,\dots$. Since $\alpha$
induces an isomorphism between the restrictions, we obtain an induced
isomorphism $\beta$ between $R_1$ and $R_2$.  It is well-known
that this means that $r = r^\prime$ where $r = r_1r_2 \dots$ and $r^\prime =
r^\prime_1 r^\prime _2 \dots$ are generalised integers.  (See [P2] for
example).
Thus we obtain the necessary condition that the integer $n$ is a divisor of the
generalised integer $r$.

We shall now improve on this necessary condition.

The isomorphism between $R|\overline{{\cal O}(x_\ast)}$ and $R|\overline{{\cal
O}(\alpha(x_\ast))}$ is given explicity by
\[
\alpha : (\dots 1,\hat{1}, x_1,x_2, \dots ) \to (\dots 1, w^\prime, y_t,
y_{t+1},\dots)
\]
where
\begin{eqnarray}
\frac{\|w^\prime\|-1}{n} + \sum^\infty_{k=1} \frac{(y_{t+k-1}-1)}{nm_{t+k-1}}
m_{t-1} = \sum^\infty_{k=1} \frac{x_k - 1}{m_k},
\end{eqnarray}
where $\|w^\prime \|$ is the cardinality of the set of points in the order
interval from the $(2t-2)$-tuple $(1,1,\dots,1)$ to $w^\prime$, and where $m_k
=
r_1r_2 \dots r_k$ for $k = 1,2\dots$.  The identity (1) follows from the fact
that there are unique canonical $R$-invariant probability measures on
$\overline{{\cal O}(x_\ast)}$ and on $\overline{{\cal O}(\alpha(x_\ast))}$ and
the quantities in (1) are the measures of the subsets $\overline{{\cal
O}(\alpha(x))}$
and $\overline{{\cal O}(x)}$ respectively.

To verify these facts one must recall how the topology of a   topological
binary relation  is  defined.  In the case of $R_1 = R|\overline{{\cal
O}(x_\ast)}$ fix      two words
\[
(x_1, x_2, \dots, x_\ell) \leq (x^\prime_1, x^\prime_2, \dots, x^\prime_\ell)
\]
in lexicographic order. Then the set $E$ of pairs
\[
((x_1,x_2,\dots,x_\ell,z_{\ell +1}, z_{\ell +2}, \dots),
(x^\prime_1, x^\prime_2, \dots, x^\prime_\ell,
z_{\ell +1}, z_{\ell +2}, \dots ))
\]
is, by definition, a basic open and closed subset for the topology. Notice that
for this set, the left and right coordinate projection maps, $\pi_\ell : E \to
\overline{{\cal O}(x_\ast)}$, $\pi_r : \to  \overline{{\cal O}(x_\ast)}$, are
injective. In the language of groupoids, $E$ is a $G$-set.  If $\lambda$ is a
Borel measure such that $\lambda(\pi_\ell(E)) = \lambda(\pi_r(E))$ for all
closed and open $G$-sets $E$, then $\lambda$ is said to be $R$-invariant. It is
easy to see that this requirement forces $\lambda$ to be the product measure
$\lambda_1 \times \lambda_2 \times \dots$ where $\lambda_k$ is the uniformly
distributed probability measure on $\{1, \dots, r_k\}$. (One can also bear in
mind that $R$-invariant measures are also $\tilde{R}$-invariant, where
$\tilde{R}$ is the topological equivalence relation (i.e. groupoid) generated
by
$R$, and that the $\tilde{R}$-invariant measures correspond to traces on the
C*-algebra of $\tilde{R}$. In our context $C^\ast(\tilde{R})$ is UHF, and the
$R$-invariant measure corresponds to the unique trace.)

  Let $\nu (x)$ denote the right hand
quantity of (1). Then the coordinates for $\alpha(x)$ are calculated from the
identity (1), bearing in mind that the ambiguity arising from the equality
$\nu(x) = \nu (x^+)$, for a gap pair $x,x^+$, is resolved by the known
correspondence of left and right gap points.

Note that if $x$ is in $\overline{{\cal O}(x_\ast)}$, and $\alpha(x) = y =
(y_k)$,
and $\| w^\prime \| = 1$ (so that $y_{-t+1}, y_{-t}, \dots, y_t$ are all equal
to 1), then, by (1),
\[
\nu(\alpha (x)) = \sum^\infty _{k=1} \frac{y_k -1}{m_k} = \sum^\infty_{k=1}
\frac{y_{t+k-1}-1}{m_{t+k-1}} = \frac{n\nu (x)}{m_{t-1}}.
\]
We have obtained the identity $\nu (\alpha (x)) = c \nu (x)$, with $c = n /
m_{t-1}$, for all points $x$ in $\overline{{\cal O}(x_\ast)}$ for which
$\nu(x)$ is small.   In fact, because of the $R$-invariance of the measures on
$\overline{{\cal O}(x_\ast)}$ and  $\overline{{\cal O}(\alpha(x_\ast))}$, which
we shall call $\lambda_1$ and $\lambda_2$ respectively, it follows that
$\nu(\alpha(x)) = c\nu(x)$ for all points $x$ for which $\alpha(x) \in
\overline{{\cal O}(x_\ast)}$. To be more precise about this, consider the left
gap points
\begin{eqnarray*}
g & = & (\dots 1, \hat{1}, 1, \dots, 1, r_{\ell + 1}, \dots ), \\
x & = & (\dots 1, \hat{1}, w, r_\ell, r_{\ell + 1}, \dots ) , \\
x^\prime & = & (\dots 1, \hat{1}, w, r_\ell - 1, r_{\ell + 1}, \dots ),
\end{eqnarray*}
where $w$ is some word $w_1, w_2, \dots, w_{\ell - 1}$. Note that the set
\[
E = \{((\dots 1, \hat{1}, w,r_\ell, z_{\ell + 1}, z_{\ell +2},\dots),
(\dots 1, \hat{1}, \dots, 1, z_{\ell + 1}, z_{\ell + 2}, \dots)) : z_j \leq
r_j\} \]
has $\pi_\ell(E) = \overline{{\cal O}(x)} \setminus \overline{{\cal
O}(x^\prime)}$ and $\pi_r(E) = \overline{{\cal O}(g)}$, and so $\nu(g) = \nu(x)
- \nu(x^\prime)$. Since $\alpha$ preserves orbits and $G$-sets we also deduce
that
\begin{eqnarray*}
\nu(\alpha(g)) & = & \lambda_1 (\overline{{\cal O}(\alpha(g))}) =
\lambda_1 (\pi_r ((\alpha \times \alpha)(E)))\\
& = & \lambda_1(\pi_\ell((\alpha \times \alpha)(E)))
= \lambda_1(\overline{{\cal O}(\alpha(x))} \setminus \overline{{\cal
O}(\alpha(x^\prime))}) \\
& = & \nu(\alpha(x)) - \nu(\alpha(x^\prime)).
\end{eqnarray*}
Thus, if we choose $\ell$ large, so that we know that $\nu(\alpha(g)) =
c\nu(g)$, we deduce that
\[
\nu(\alpha(x)) - \nu(\alpha(x^\prime))
= \nu(\alpha(g)) =
c \nu(g) =
c(\nu(x) - \nu(x^\prime)),
\]
from which it follows that $\nu(\alpha(x)) = c(\nu (x))$ for general points $x$
with $\alpha(x)$ in $\overline{{\cal O}(x^\prime)}$.

  We can similarly
extend this identity to points in the set
\[
X_0 = \{(y_k) \in X : \exists k_0 \mbox { such that } y_k = 1 \mbox { for all }
k \leq k_0\}
\]
and the extension of $\nu$ given by
\[
\nu (y) = \sum^\infty_{k=1} (y_{-k} - 1)s_0s_1\dots s_{k-1} + \sum^\infty_{k=1}
\frac{y_k-1}{m_k}
\]
for $y$ in $X_0$, where $s_0 = 1$.   The range of $\nu$ on the gap points of
$X_0$ is the
additive cone  of rationals of the form $\ell /m_k$ for some $k = 1,2,\dots$
and
some natural number $\ell$. The identity $\nu (\alpha (x)) = c \nu (x)$ for $x$
in $X_0$ shows that multiplication by $c$ is a bijection of the cone. From this
we obtain the necessary condition that $c$ has the form
\[
c = p^{a_1}_1 \dots p ^{a_d}_d
\]
where $a_i \in \ZZ, 1 \leq i \leq d$, and where $p _1, \dots p _d$ are
primes which divide infinitely many terms of the sequence $(r_k)$.

We now improve further on this condition by considering the fact that $\alpha$
is a homeomorphism of $X$ and is determined by its restriction to $X_0$.

Suppose, by way of contradiction, that $a_1 \neq 0$ and that $p_1$ does not
divide infinitely many terms of the sequence $(s_k)$.
Note that $c$ only depends on $\alpha$, thus, replacing $\alpha$ by its
inverse if necessary, we may assume that $a_1 > 0$.
By relabelling we may also
assume that $p_1$ divides no terms of the sequence.  Without loss of
generality assume that $s_1 > 1$ and consider the proper clopen subset $E$ of
points $y = (y_k)$ in $X$ with $y_{-1} = 1$.
We show that $\alpha (E)$ is dense, which is the
desired contradiction.	Observe first that the range of $\nu$ on $E \cap X_0$ is
the union of the intervals $[ks_1, ks_1 + 1]$ for $k = 0,1,2,\dots$ Pick $x$ in
$X_0$ arbitrarily, pick $j$ large, and consider the countable set
\[
F_j(x) = \{x^\prime \in X_0 : x^\prime = (x^\prime _k) \mbox { and } x^\prime
_k
= x_k \mbox { for all } k \geq -j\}.
\]
The range of $\nu$ on $F_j(x)$ is an arithmetic progression of period
$s_1 s_2 \dots s_j$. In view of the identity $\nu(\alpha (y)) = c \nu (y)$, the
range of $\nu$ on $\alpha (E) \cap X_0$ is the union of the intervals $[cks_1,
cks_1 + c]$, which is an arithmetic progression of intervals of period $cs_1$.
It follows from our hypothesis on $p_1$ that one of these intervals contains
a point in $\nu(F_j(x))$, and so $\alpha (E)$ meets $F_j(x)$.
 Since the intersection of the sets $F_1(x), F_2(x), \dots$
is the singleton $x$, it follows that $x$ lies in the closure of $\alpha(E)$.
Since $X_0$ is dense it follows that $\alpha(E)$ is dense as desired.

We have now shown that if $\alpha$ is an automorphism of $R = R(r_k,s_k)$, then
$\nu(\alpha (x)) = c\nu (x)$ for all $x$ in $X_0$ where $c$ has the form $c =
p^{a_1}_1 p^{a_2}_2 \dots p^{a_d}_d$ where $a_1, \dots, a_d$ are
integers and where $p_1, \dots p _d$ are primes which divide infinitely
many terms of $(r_k)$ and of $(s_k)$.  It is also clear from the above that for
each such $c$
there is at most one automorphism $\alpha$ satisfying the identity $\nu (\alpha
(x)) = c \nu (x)$.
It follows that the map
\[
\alpha \to (a_1, \dots, a_d)
\]
is an injective group homomorphism from Aut$R$ to $\ZZ ^d$. ($d$ may be
infinite.)  It remains to show that this map is surjective.  One way to do this
is to start with $c$ of the required form above and to show that the bijection
of $X_0$ induced by multiplication by $c$ (that is, the bijection $\alpha$
satisfying $\nu (\alpha (x) ) = c \nu (x) )$ does extend to an order preserving
homeomorphism of $X$ which defines an automorphism of $R$.  Another way, which
we now follow, is to make the connection between $R(r_k, s_k)$ and
$T(r_k,s_k)$,
and to determine generators of Aut$R$ in terms of commuting diagrams, as we
indicated after  the statement of Theorem 1.

Consider the diagram
\[
\begin{array}{llllllllll}
\IC &
\stackrel{\rho_{r_1}}{\to} & M_{r_1}
& \stackrel{\sigma_{s_1}}{\to}
& M_{s_1} \otimes M_{r_1} &
\stackrel{\rho_{r_2}}{\to} & M_{s_1} \otimes M_{r_1} \otimes M_{r_2} &
\stackrel{\sigma_{s_2}}{\to} & \dots & B\\
\uparrow & & \uparrow & & \uparrow & & \uparrow & & & \uparrow\\

\IC &
\stackrel{\rho_{r_1}}{\to} & T_{r_1}
& \stackrel{\sigma_{s_1}}{\to}
& T_{s_1r_1}              &
\stackrel{\rho_{r_2}}{\to} & T_{s_1r_1r_2}                        &
\stackrel{\sigma_{s_2}}{\to} & \dots & A\\
\end{array}
\]

The vertical maps are inclusions, where $T_{s_1r_1r_2}$, for example, is
realised in terms of the lexicographic order on the indices $(i,j,k)$ of the
minimal projections $e_{ii} \otimes e_{jj} \otimes e_{kk}$ in $M_{s_1} \otimes
M_{r_1} \otimes M_{r_2}$. (For more detail concerning this discussion, read the
introduction of [HP].) The maximal ideal space of the diagonal C*-algebra $A
\cap A^\ast$ is naturally identified with the space $X$.  Indeed, $x = (x_k )$
in $X$ corresponds to the point in the intersection of the Gelfand supports of
the projections
\[
e(x,N) = e_{x_{-N,-N}} \otimes \dots \otimes e_{x_{-1,-1}} \otimes e_{x_{1,1}}
\otimes \dots \otimes e_{x_{N,N}}
\]
for $N = 1,2,\dots$. Furthermore, $(x,y)$ belongs to $R = R(r_k,s_k)$ if and
only if for all large $N$ there is a matrix unit in the appropriate upper
triangular matrix algebra with initial projection $e(y,N)$ and final projection
$e(x,N)$. (In fact $R$ is the fundamental relation of the limit algebra $A$.)

Suppose now that $r_k = s_k = p$ for all odd $k$ and let $\alpha$ be the
isometric  automorphism of $T(r_k, s_k)$ determined by the diagram given in the
introduction.  Let $\alpha$ also denote the induced automorphism of $R$.  We
prove that $\nu (\alpha (x)) = p ^{-1} \nu (x)$, completing the proof of the
theorem.

Let us calculate $\alpha (e (x, N)) $, where $N$ is even, $x = (\dots, 1,
\hat{1}, 2, 1, \dots)$, and where we abuse notation somewhat and write $e(x
, N)$
for the image of $e(x	  , N)$ in the limit algebra.  Let $d(N) = s_N \dots s_1
r_1 \dots r_N$, and let $e(x	 , N)$ occupy position $a(N)$ in the
lexicographic ordering of the $d(N)$ matrix units.  Consider the following part
of the diagram defining $\alpha$.

\[
\begin{array}{lllll}
T_{s_N..r_N} & \stackrel{\rho_{r_p}}{\to}
& T_{s_N..r_{N+1}} & \stackrel{i}{\to} & A\\
& \sigma_p & & & \downarrow\\
& & T_{s_N..r_{N+1}}  & \stackrel{i}{\to} & A
\end{array}
\]

Then
\[
\rho_p (e(x,N)) = \sum^p _{k=1} e(x,N) \otimes e_{kk}.
\]
On the other hand $\sigma_p (e (x,N))$ is the summation of the diagonal matrix
units in positions $a(N), a(N) + d(N), \dots, a(N) + (p-1) d(N)$ in the
lexicographic order.   Let these projections correspond to the matrix unit
tensors with subscripts $z^{(i)} = (z^{(i)}_{-N},\dots z^{(i)}_{N+1})$ for $1
\leq i \leq p$, and denote the projections themselves by $f_1, \dots f_p$,
respectively.  It follows (from the partial diagram above) that the
homeomorphism $\alpha : X \to X$ maps the support of $e(x,N)$ onto the union of
the supports of $f_1, \dots, f_p$.  Denote these supports by $E(x,N), F_1,
\dots, F_p$ respectively.  Since $X_0$ is invariant for $\alpha$,
\[
\alpha (E(x,N) \cap X_0) = \bigcup^p_{k=1} F_k \cap X_0.
\]
Notice that $x$ is the unique point in $E(x,N) \cap X_0$ with the property that
if $y \in E(x,N) \cap X_0$ and $(x,y) \in \tilde{R}$ then $(x,y) \in R$.  The
point in the union of $F_1 \cap X_0, \dots, F_p \cap X_0$ with this minimum
property is the point
\[
u = (\dots 1 \ 1 \ z^{(1)}_{-N} , \dots, z^{(1)}_{N+1}, 1, 1, \dots )
\]
and so $\alpha (x) = u$.  Finally one can verify that $\nu (x) = p ^{-1}$ and
$\nu (u) = p ^{-2}$, as desired.
\end{mproof}

Recall that the fundamental relation $R(A)$ of a canonical triangular
subalgebra
$A$ of an AF C*-algebra $B$ is the topological binary relation on the Gelfand
space $M(A \cap A^\ast)$ induced by the partial isometries of $A$ which
normalise $A \cap A^\ast$. (See [P2].) In [HP] we identified $R(A)$, for $A =
T(r_k,s_k)$,
with $R(r_k,s_k)$.  (This identification is also effected in the proof above
by virtue of the fact that a matrix unit system determines $R(A)$.) Let $\beta$
be an isometric automorphism of $A$. Then
$\beta$ induces an automorphism of $R(A)$ (because $\beta(A \cap A^\ast) = A
\cap A^\ast$ and $\beta$ maps the normaliser onto itself). Thus $\beta$
determines an automorphism of $R(r_k,s_k)$ and so by the last theorem there is
an isometric automorphism $\alpha$ of $A$ such that $\gamma = \alpha^{-1} \circ
\beta$ induces the trivial automorphism of $R(r_k,s_k)$. This means that
$\gamma$ is an isometric automorphism with $\gamma$ equal to the identity map
on
$A \cap A^\ast$.

{\bf Lemma}
Let $\gamma$ be an automorphism of $T(r_k,s_k)$ which is the identity on the
diagonal subalgebra (and which is not necessarily isometric). Then $\gamma$ is
approximately inner.

\begin{mproof} Let $A = T(r_k,s_k)$ and let $A_1 \to A_2 \to \dots $ be the
direct system defining $A$.  The hypothesis is that $\gamma(c) = c$ for all $c$
in $C = A \cap A^\ast$. This ensures that $\gamma (\tilde{A}_n) = \tilde{A}_n$
where $\tilde{A}_n$ is the subalgebra generated by $A_n$ and $C$. To see
this,recall from Lemma 1.2 of [P1] that there are contractive  maps $P_n : A
\to
\tilde{A}_n$ which are defined in terms of limits of sums of compressions by
projections in $C$, and so, for $a$ in $\tilde{A}_n$, $\gamma (a) =
\gamma(P_n(a)) = P_n (\gamma (a))$. The restriction automorphism $\gamma|
\tilde
{A}_n$ is necessarily inner. Indeed identify $\tilde{A}_n$ with $T_r \otimes
D$,
for appropriate $r$, where $D$ is an abelian approximately finite C*-algebra
and
let $u_i \in D$, $1 \leq i \leq r-1$, be the invertible elements such that
$\gamma(e_{i,i+1}) = e_{i,i+1} \otimes u_i$. Also set $u_0 = 1$.  Then it
follows that $\gamma (a) = u^{-1} au$, where
\[
u = \sum^r_{i=1} e_{i,i} \otimes u_0u_1\dots u_{r-1}
\]
Furthermore, since $\gamma (e_{1,r}) = e_{1,r} \otimes u_0u_1\dots u_{r-1}$, it
follows that $\|u\| \leq \| \gamma \|$. Similarly $\|u^{-1}\| \leq \|
\gamma^{-1}\|$.  The inner automorphisms $Adu^{-1}$, for varying $n$, thus form
a uniformly bounded sequence which converge pointwise on each $A_n$, and so
determine an approximately innder automorphism.
\end{mproof}

It follows from Lemma 1 and the preceeding discussion that
\[
Aut_{isom} A/I(A) = Aut R(A) = \ZZ ^d .
\]

{\bf Remark 1.}
Suppose that $\delta \in Aut$$A$. Then $\delta$ determines a scaled group
homomorphism $\delta_\ast : K_0 (A) \to K_0(A)$ which preserves the algebraic
order on the scale $\Sigma(A)$ of $K_0(A)$. Thus, by the main theorem of [P3],
(which can also be found in [P4])
there is an isometric algebra automorphism of $A$, $\phi$ say, with $\phi_\ast
=
\delta_\ast$.  In particular $\psi = \phi ^{-1} \circ \delta$ has $\psi_\ast$
trivial. This means that if $P : A \to A \cap A^\ast$ is the diagonal
expectation, then $P(\psi(e)) = e$ for each projection $e$ in $A \cap A^\ast$.
Thus to show that $AutA/I(A) = \ZZ ^d$ it remains only to show that such
automorphisms $\psi$ are approximately inner.

{\bf Remark 2.}
There are approximately inner automorphisms of alternation algebras which are
not inner. To see this, consider the standard limit  algebra $A =
{\displaystyle\lim_\to}(T
_{2^n}, \sigma)$.

Let $\lambda$ be a unimodular complex number and let
$d_n = \lambda e_{1,1} + \lambda^2e_{2,2} + \dots + \lambda^{2^n}e_{2^n,2^n}$.
Then $d_nad_n^{-1} = d_mad_m^{-1}$ if $a \in T_{2^n}$
and $m > n$, from
which it follows that $\alpha(a) = \lim_n (d_n a d^{-1}_n)$ is
an isometric approximately inner automorphism.

Suppose now that $\alpha$ is inner, and $\alpha(a) = gag^{-1}$ for some
invertible $g$ in $A$. Since $\alpha(c) = c$ for all $c$ in the masa $C$ it
follows that $g \in C$.  In particular $\|\alpha - \beta \| \leq \frac{1}{4}$
for some inner automorphism $\beta$ of the form $\beta(a) = h ah^{-1}$ where,
for some large enough $n$,
$h \in T_{2^n} \cap (T_{2^n})^\ast$.
However, in $T_{2^m}$, for large $m$, the
diagonal element $h$ has matrix entries
which are periodic with period $ 2^{n}$. One can now verify that if
$\lambda$ is chosen so that no power of order $2^k$ is unity then for
large enough $m$ there exist matrix units $e \in T_{2^m}$ such that
$\| \lambda e - heh^{-1}\| > \frac{1}{4}$, a contradiction.


{\bf Remark 3.}
Let $(x,y)$ be a point in $R(C^\ast(A(r_k,s_k)))$ with $x = (\dots, x_{-2},
x_{-1}, x_1, x_2, \dots)$, $y = (\dots, y_{-2}, y_{-1}, y_1, y_2, \dots)$.
Then,
although $\nu(x)$ and $\nu(y)$ may be infinite, we may define $d(x,y)$ as the
sum
\[
\sum^\infty_{k=1} (y_{-k} - x_{-k})s_0 s_1\dots s_{k-1} +
\sum^\infty_{k=1} \frac{y_k - x_k}{r_1r_2\dots r_k}
\]
because only finitely many terms are nonzero.  Since $d(x,y) = d(x,z) +
d(z,y)$,
and $(x,y) \in R(r_k,s_k)$ if and only if $d(x,y) \geq 0$, it follows that
$d(x,y)$ is a continuous real valued cocyle determining $A(r_k,s_k)$ as an
analytic subalgebra of $C^\ast(A(r_k, s_k))$.  See [V], where some special
cases
are discussed as well as some general aspects of analyticity.

\it  Added Dec 1992 : Unfortunately the proof of the classification
of alternation algebras given in [HP] and [P4] appears to be incomplete.
(It is not clear, in [P4], whether $q$ can be chosen with the desired
properties.) However the present paper is independent of [HP] and
the arithmetic progression argument above can be adapted, to the case
of an isomorphism $\alpha$ between two alternation algebras, to show
that the supernatural numbers for the standard multiplicities are
finitely equivalent.

\rm
\vspace{1cm}

{\bf References}

[HP] \quad A. Hopenwasser and S.C. Power, Classification of limits of
triangular
matrix algebras, Proc. Edinburgh Math. Soc., to appear.

[P1] \quad S.C. Power, On ideals of nest subalgebras of C*-algebras, Proc.
London Math. Soc., 50 (1985), 314-332.

[P2] \quad S.C. Power, The classification of triangular subalgebras of AF
C*-algebras, Bull. London Math. Soc., 22 (1990), 269-272.

\newpage

[P3] \quad S.C. Power, Algebraic orders on $K_0$ and approximately finite
operator algebras, preprint 1989, to appear in J. Operator Th.

[P4] \quad S.C. Power, Limit Algebras : An Introduction to Subalgebras
of C*-algebras, Pitman Research Notes in Mathematic Series,
No. 278, Longman, 1992.

[Po] \quad Y.T. Poon, A complete isomorphism invariant for a class of
triangular
UHF algebras, preprint 1990.

[V] \quad B.A. Ventura, Strongly maximal triangular $AF$ algebras,
International J. Math., 2 (1991), 567--598.

 \end{document}